\documentclass[11pt]{article}
\usepackage{amsmath,amsthm,amsfonts,amscd,bezier}
\usepackage{amssymb}
\usepackage[usenames]{color}
\newtheorem{theorem}{Theorem}%[section]
\newtheorem{lemma}[theorem]{Lemma}

\newtheorem{corollary}[theorem]{Corollary}

\newtheorem{proposition}[theorem]{Proposition}
\newcommand{\pf}{\noindent{\bf Proof:\ \ }}
\newcommand{\cqd}{{\hfill $\rule{2mm}{2mm}$}\vspace{1cm}}

\def\C{\mathbb{C}}

\def\B{{\cal B}}

\def\A{{\cal A}}

\def\H{{\cal H}}

\begin{document}

\pagestyle{plain} \pagenumbering{arabic}
\title{The Analytic Classification of Plane Curves with Two Branches}
\author{Hefez, A., Hernandes, M.E. and Rodrigues Hernandes, M.E.
\thanks{The first two authors
were partially supported by CNPq grants}}
\date{ \ } 
\maketitle

\begin{center} Abstract
\end{center}

{\small In this paper we solve the problem of analytic classification of plane curves singularities with two branches by presenting their normal forms. This is accomplished by means of a new analytic invariant that relates vectors in the tangent space to the orbits under analytic equivalence in a given equisingularity class to K\"ahler differentials on the curve.}

\section{Introduction}

Let $(f):f=0$ be the germ of a reduced plane analytic curve, that is, the curve associated to a reduced element $f$ in $\C\{X,Y\}$, the ring of convergent power series in two variables over the complex numbers. Mather's contact equivalence asserts that $f$ and $g$ are equivalent, writing $(f)\sim (g)$, if and only if there exist $\Phi\in Aut(\C\{X,Y\})$ and a unit $u$ in $\C\{X,Y\}$ such that $\Phi(f)=ug$.

The aim of this work is to initiate the analytic classification of germs of reducible (but reduced) plane curves, that is, the classification for Mather's contact equivalence. The irreducible case was solved by the first two authors in \cite{HH} and our results here concern curves with two components.

From now on, we will assume that $f$ has two irreducible components $f_1$ and $f_2$. Each branch $(f_i)$ admits a parametrization $\phi_i\colon (\C,0) \to (\C^2,0)$. We will use coordinates $t_1$ and $t_2$ on $(\C,0)$ (one for each $\phi_i$) and coordinates $x,y$ on $(\C^2,0)$ (the same for both). Now, because each branch is invariant by changes of coordinates in the source of the $\phi_i$, and the curve is analytically invariant by any automorphim of $(\C^2,0)$ (the same automorphism for both branches), we easily conclude that contact equivalence for curves $(f)$ is translated into $\A$-equivalence on the associated bigerms $\phi =[\phi_1, \phi_2]$, i.e., changes of analytic coordinates in the source and in the target. The space of bigerms will be denoted by $\B$.

The $\A$-equivalence in $\B$ is induced by the action of the group $\A =  Aut(\C\{t_1\}) \times  Aut(\C\{t_2\}) \times Aut(\C\{X,Y\})$,
as follows: 
$$
(\rho_1,\rho_2,\sigma)\cdot \phi= [\sigma \circ \phi_1 \circ \rho_1^{-1},  \sigma \circ \phi_2
\circ \rho_2^{-1}].$$

Our analysis will be splitted into two cases, namely, whether the two components of $(f)$ have distinct tangents (the transversal case) or equal tangents. In what follows, we will denote by $m_i$ the multiplicity of $f_i$, $i=1,2$.\smallskip

\noindent Case 1) Distinct tangents. In this case, by $\A$-equivalence, we may assume that the tangent of the first component is $(Y)$ and of the second one is $(X)$, so that 
\[\phi_i=(x(t_i),y(t_i)), \ \hbox{where} \ {\rm ord}_{t_1}x(t_1)<{\rm ord}_{t_1}y(t_1) \ \hbox{and} \
{\rm ord}_{t_2}x(t_2)>{\rm ord}_{t_2}y(t_2).
\] 

\noindent Case 2) Same tangent. In this case, by $\A$-equivalence, we may assume that the common tangent is $(Y)$, in which case, $\phi_i=(x(t_i),y(t_i))$ with ${\rm ord}_{t_i}x(t_i)<{\rm ord}_{t_i}y(t_i)$, $i=1,2$.\smallskip

To describe the elements of $\A$ that preserve the tangent cone of the bigerm, it is convenient to introduce the subgroup $\H$  of $\A$ of homotheties:
$$\H= \{ (\rho_1,\rho_2,\sigma)\in \A; \ j^1 \rho_i=\alpha_it_i, \ j^1\sigma=(ax,by), \  \alpha_i, a, b \in \C^*, \ i=1,2\},
$$
where $j^k \xi$ is the $k$-th jet of any $n$-tuple $\xi$ of power series.

Now, the elements of $\A$ we are looking for are the compositions $h\circ g$, where $h\in \H$ and $g$ belongs to the classical group $\A_1$, if we are in Case (1) or $g$ belongs to $\tilde{\A_1}$, if we are in Case (2), where 
\[
\tilde{\A_1}=\{(\rho_1,\rho_2, \sigma)\in \A; \ j^1 \rho_i=t_i,\ i=1,2, \ j^1 \sigma=(x+by,y), \ b\in \C
\}.
\]

Notice that the group $\A_1$ is the subgroup of elements of $\tilde{\A_1}$ with $b=0$.

The strategy we use for our classification is to first analyze the action of $\A_1$, or $\tilde{\A_1}$, on the elements of $\B$, according they belong, respectively, to Case (1) or to Case (2) and then to take into account the homotheties. 

To find distinguished representatives in each case, under the action of the corresponding group, we will use the Complete Transversal Theorem (cf. \cite{BKP}).\smallskip

\noindent  {\sc The Complete Transversal Theorem} ({\em CTT}). \ {\em Let $G$ be a Lie group acting on an affine space $\mathbb A$ with underlying vector space $V$ and let $W$ be a subspace of $V$. Suppose that $v \in V$ is such that $TG(v+w)=TG(v)$, $\forall\; w\in W$, where the notation $TG(z)$ means the tangent space at $z$ of the orbit $G(z)$, as vector subspace of $V$. If $W\subset TG(v)$, then  $G(v+w)=G(v)$, $\forall\; w\in W$.} \smallskip

We denote by $\B^k$ the vector space of $k$-jets of elements of $\B$ and by $G^k$ the Lie group of $k$-jets of elements of $G$, where $G$ is one of the groups $\A_1$ or $\tilde{\A_1}$.

We will show, in the next proposition, that the hypothesis of {\em CTT} holds for an element $j^k\phi \in \B^k$ which is a bigerm as in Case (2), where $\phi\in \B$,  and for $W=H_\phi^k$, the subspace of homogeneous elements of degree $k$ of $\B^k$ such that the two components, as a bigerm, of the elements in $j^k\phi+H_\phi^k$ have all same multiplicity and same tangent, that is, 
$$H_\phi^k=\left\{[a_1t_1^k,b_1t_1^k), (a_2t_2^k,b_2t_2^k)] \in \B^k; \ a_i=b_i=0, \ \text{if} \ k\leqslant m_i, \ i=1,2\right\}.$$

We describe below the elements of the tangent spaces to the orbit of $j^k\phi$ in $\B^k$ under the actions of the groups $\A_1^k$  and $\tilde{\A}_1^k$:
\begin{equation}\label{tangentvector}
j^k\left[ (\phi_{11}'\epsilon_1+\eta_1(\phi_1), \phi_{12}'\epsilon_1+\eta_2(\phi_1)),\; 
(\phi_{21}'\epsilon_2+\eta_1(\phi_2), \phi_{22}'\epsilon_2+\eta_2(\phi_2) )
\right],
\end{equation}
where $\phi_i=\big( \phi_{i1}, \phi_{i2}\big)$ the $(')$ sign means derivative with respect to the corresponding parameter, $\epsilon_i\in (t_i)^2\C\{t_i\}$, $i=1,2$, $\eta_2\in (x,y)^2\C\{x,y\}$ and 

\noindent a) $\eta_1\in (x,y)^2\C\{x,y\}$, in the $\A_1^k$ case, or

\noindent b) $\eta_1\in (x^2,y)\C\{x,y\}$, in the $\tilde{\A}_1^k$ case. 

The case of $\A_1^k$ is classically known (cf. \cite{Gi}) and the other one can be computed in a similar way.

\begin{lemma} If $\phi\in \B^k$ as in case (2), $(\rho_1,\rho_2,\sigma)\in \tilde{\A}_1^k$, with $j^1\sigma=(x+by,y)$,  and $\psi\in H_\phi^k$, then
\[
j^k[(\rho_1,\rho_2,\sigma)\cdot (\phi+\psi)]=j^k[(\rho_1,\rho_2,\sigma)\cdot \phi]+\psi + \theta,
\]
where $\theta=[(bc_1t_1^k,0),(bc_2t_2^k,0)]$, with $b,c_1,c_2\in \C$, depending only upon $\psi$.
\end{lemma}
The proof is straightforward, following easily from the definitions.

\begin{proposition} If $\phi\in B^k$ as in case (2) and $\psi\in H_\phi^k$, then
\[
T\tilde{\A}_1^k(\phi+\psi)=T\tilde{\A}_1^k(\phi).
\]
\end{proposition}
\pf Recall that $T\tilde{\A}_1^k(\psi)$ of an element $\psi\in \B^k$ is given by the image of the differential at the identity $I$ of the map $\ \Phi_{\psi}\colon \tilde{\A}_1^k\to \tilde{\A}_1^k (\psi)$, $\Phi_{\psi}(g)=g\cdot \psi$.

Therefore, any vector in $T\tilde{\A}_1^k(\psi)$ is of the form $(\Phi_{\psi}\circ
\lambda)'(0)$, where $\lambda \colon (-\alpha, \alpha) \to \tilde{\A}_1^k$, $\lambda(u)=(\rho_{1u},\rho_{2u},\sigma_u)$ is a curve in $\tilde{\A}_1^k$ such that $\lambda(0)=I$. Notice that since $\lambda(u)\in \tilde{\A}_1^k$, then $j^1\sigma_u= (x+b(u)y,y)$.

As a consequence of the above discussion, and from the previous lemma, we have that
\[ \begin{array}{rcl}
 (\Phi_{\phi + \psi} \circ \lambda)'(0)&=&\underset{u \to 0}{lim}\
 \displaystyle\frac{\lambda(u)\cdot(\phi + \psi)-\lambda(0)\cdot (\phi + \psi)}{u} =\underset{u \to 0}{lim}\
 \displaystyle\frac{\lambda(u)\cdot \phi+\psi+\theta(u) - \phi - \psi}{u} \vspace{.3cm}\\
&=&  (\Phi_{\phi}\circ \lambda)'(0) \ + \ \underset{u \to 0}{lim}\
 \displaystyle\frac{\theta(u)}{u}=(\Phi_{\phi}\circ \lambda)'(0) \ + \  \theta'(0),  
 \end{array}
 \]
where $\theta(u) = \left[(b(u) c_1t_1^k ,0), (b(u) c_2t_2^k ,0)\right]$, with $b(0)=0$ and $\theta(0)=0$, since $\lambda(0)=I$.

Taking in the description of the tangent spaces to the orbits in $\B^k$ under the $\tilde{\A}_1^k$-action, $\eta_i=0$ and $\epsilon_i(t_i)=\frac{b'(0)c_i}{m_{i}}t_i^{k-m_i+1}$, $i=1,2$, one may easily check that $\theta'(0)\in T\tilde{\A}_1^k(\phi+\psi)\cap T\tilde{\A}_1^k(\phi)$, $\forall\; \psi \in H_\phi^k$.   
\cqd

Notice that our proof may be, without any extra effort, extended to multigerms, and contains as an immediate corollary the result for the $\A_1^k$-action (just take $b=0$). Remark also that the result may be used to make substantial simplifications in the arguments in Section 5 of \cite{HH}.

At this point it will be convenient to unify the notation for both actions $\A_1$ and $\tilde{\A_1}$.
We define $\A[1]=\A_1$ and $\A[2]=\tilde{\A_1}$, which we condense in the notation $\A[\delta]$, $\delta=1,2$. Notice that if $\phi$ and $\varphi$ are $\A^k[\delta]$ equivalent, then $H_\phi^k=H_\varphi^k$. Observe also that if $\phi\in \B$ is in Case ($\delta$) ($\delta=1,2$), then  $\phi+\psi$ is also in Case ($\delta$), for all $\psi\in H_\phi^k$.

\section{Normal Forms}

Given an element $\phi\in \B$, we are looking for elements $\psi\in H_\phi^k$ such that $j^k(\phi+\psi)=j^{k-1}\phi$ and $\phi+\psi$ is $\A^k[\delta]$ equivalent to $\phi$. So that in this way we will be able to eliminate terms of order $k$ in $\phi$ without changing neither its $k-1$ jet nor its equivalence class. From the {\em CTT} it is sufficient to verify when an element $\psi \in H_\phi^k$ belongs to the tangent space to the orbit of $\phi$ under the action of the group $\A^k[\delta]$. Similarly, as in the proof of Proposition 2 we get that 
\[
[(at_1^k,0),((\delta-1)bt_2^k, (2-\delta)bt_2^k)] \in H_\phi^k\cap T\A^k[\delta](\phi), \ a,b\in \C.
\]

With the above considerations, we have that any bigerm is $\A$-equivalent to a bigerm $\phi=[\phi_1,\phi_2]$ in Puiseux form, that is, $\phi_1=(t_1^{m_1}, \sum_{i>m_1}a_{1i}t_1^i)$ and \smallskip

\noindent Case 1) Distinct tangents: $\phi_2=(\sum_{i>m_2}a_{2i}t_2^i, t_2^{m_2})$;\smallskip

\noindent Case 2) Same tangent: $\phi_2=(t_2^{m_2},\sum_{i>m_2}a_{2i}t_2^i)$.\smallskip

The pair $m=(m_1,m_2)$ will be referred to as the multiplicity of the bigerm $\phi$.

In order to get more refined parametrizations for a bigerm we have to impose some restriction on it. This is done by fixing analytic invariants. 

As a first invariant we consider the semigroup of values
$$\Gamma = \left\{\nu(\eta):=\left(\nu_1(\eta),\nu_2(\eta)\right); \ \eta \in \C\{x,y\}
  \right\},$$
where $\nu_i(\eta)=ord_{t_i}(\eta \circ \phi_i)$, $i=1,2$. This invariant characterizes completely the topological type of the curve as an immersed germ at the origin of the plane (cf. \cite{W} or \cite{Ga}). Two curves having same $\Gamma$ invariant are called equisingular. 

Fixing the semigroup of values, which determines the intersection index of the two  branches of the curve, we are fixing the contact order of their parametrizations. This will imply the coincidence of the coefficients of the Puiseux expansions of the branches up to the order of contact minus 1. On the other hand, since $\Gamma$ has a conductor $(c_1,c_2)$, we may eliminate analytically all terms in both parametrizations with order greater than $c-1$, where $c=\max\{c_1,c_2\}$, without affecting the preceding terms (cf. \cite{Ga}). This tells us that we have simultaneous finite determinacy of both parametrizations and gives us a finite dimensional space of parameters $\Sigma_\Gamma$ for a complete set of analytic representatives in the equisingularity class determined by $\Gamma$.  

With the semigroup $\Gamma$, we get only a rough normal form for bigerms. In order to refine this normal form, we will use the finer analytic invariant
$$  
\Lambda = \left\{\nu(\omega):=\left(\nu_1(\omega),\nu_2(\omega)\right); \ \omega \in \C\{x,y\}dx+\C\{x,y\}dy
  \right\},$$
where for $\omega=\eta_1dx+\eta_2dy$ with $\eta_i\in \C\{x,y\}$, $i=1,2$, we define
\[
\nu_i(\omega) :=ord_{t_i} \omega(\phi_i)+1 =ord_{t_i}(\eta_1(\phi_i)\phi'_{i1}+\eta_2(\phi_i)\phi'_{i2})+1.
\]
The fact that $\Lambda$ is an analytic invariant is clear since by its definition it is independent from reparametrizations of the branches and change of coordinates in $\C^2$. From the definition it also follows that $\Gamma\setminus \{(0,0)\}\subset \Lambda$.

It is easy to check that the set $\Lambda$ has the following properties: 

\noindent A) If $(a_1,a_2),(b_1,b_2)\in \Lambda$, are such that $a_1<b_1$ and $a_2>b_2$, then $(a_1,b_2)\in \Lambda$.

\noindent B) If $(a_1,a_2),(a_1,b_2)\in \Lambda$, then there exists $(a, \min\{a_2,b_2\})\in \Lambda$ with $a> a_1$. The same is true reversing the roles of the axes.

This is sufficient to guarantee that $\Lambda$ behaves combinatorially as $\Gamma$, except that it is not a semigroup. In $\Lambda$ there is a finite subset $M$ of points $(k_1,k_2)$, called the maximal points of $\Lambda$, contained in the rectangle with sides parallel to the axes and opposite vertexes the origin of $\mathbb N^2$ and the conductor $(c_1, c_2)$ of $\Gamma$, such that $F_1(k_1,k_2)=F_2(k_1,k_2)=\emptyset$, where for $(a_1,a_2)\in \mathbb N^2$,
$$F_i(a_1,a_2)=\{ (b_1,b_2)\in \Lambda; \, a_i=b_i, \ 
b_j>a_j, \ i\neq j\},$$
and respectively called the vertical and the horizontal fibers of $(a_1,a_2)$.

In particular, the set $\Lambda$ is determined by the 
sets of values of differentials $\Lambda_1, \Lambda_2$ of the branches of the curve and the maximal points of $\Lambda$  (cf. \cite{Ga} or \cite{D}, in the case of the set $\Gamma$). 

This implies that there are finitely many possibilities for sets of values of differentials $\Lambda$ for each equisingularity class of curves.

There is a tight connection between the tangent space to the orbit of a bigerm $\phi$ under the action of the group $\A[\delta]$ and the set 
\[
\Lambda[\delta]=\{\nu(\omega)- m; \; \omega \in \Omega[\delta]\} \subset \Lambda -m, 
\]
where $m=(m_1,m_2)$ is the multiplicity of the bigerm $\phi$ and
\[ \Omega[\delta]=\{\eta_1dx+(\beta(\delta-1)y+\eta_2)dy; \;
\eta_1,\eta_2\in (x,y)^2\mathbb{C}\{x,y\}, \beta\in\mathbb{C}\}.
\]

The same argument used for $\Lambda$ shows that the set $\Lambda[\delta]$ is an invariant with respect to the action of the group $\A[\delta]$. For each fixed semigroup of values $\Gamma$ there exist finitely many possibilities for $\Lambda[\delta]$. The finite dimensional space that parametrizes the bigerms in Puiseux form with fixed $\Gamma$ and $\Lambda[\delta]$ will be denoted by $\Sigma_{\Gamma,\Lambda[\delta]}$, which we will identify with the set of the bigerms that they determine. 

\begin{proposition}
Let $\phi = [\phi_1, \phi_2] \in
\Sigma_{\Gamma,\Lambda[\delta]}$ with $\delta =1,2$. Given $h_i\in\mathbb{C}\{t_i\},\ i=1,2$, we have that 
$[(0,h_1),(-1)^{\delta}((2-\delta)h_2,(\delta-1)h_2)]\in T\A^k[\delta] (\phi)$
 if and only if there exists $\omega \in \Omega[\delta]$ such that 
$h_i = j^k \frac{\omega \circ \phi_i}{m_it_i^{m_i-1}}$.
\end{proposition}
\pf If $[(0,h_1),(-1)^{\delta}((2-\delta)h_2,(\delta-1)h_2)]\in T\A^k[\delta](\phi)$, then from (\ref{tangentvector}) there exist
$\epsilon_i \in (t_i^2)\mathbb{C}\{t_i\}, \ i=1,2, \ \eta,\eta_2 \in
(x,y)^2\mathbb{C}\{x,y\},\ \eta_1=\beta(\delta-1)y+\eta$ with 
$\beta\in\mathbb{C}$ such that \smallskip

$0=\phi_{11}'\cdot \epsilon_1 + \eta_1(\phi_1)\ mod\ t_1^{k+1}$,\smallskip

$h_1= \phi_{12}' \cdot \epsilon_1 + \eta_2(\phi_1)\ mod\ t_1^{k+1}$,\smallskip

$(2-\delta)(-1)^{\delta}h_2=\phi_{21}'\cdot \epsilon_2 + \eta_1(\phi_2)\ mod\ t_2^{k+1}$, \ and\smallskip

$(\delta-1)(-1)^{\delta}h_2= \phi_{22}' \cdot \epsilon_2 +
\eta_{2}(\phi_2)\ mod\ t_2^{k+1},$ \smallskip

\noindent that is,
$j^k \epsilon_1=-j^k \frac{\eta_1(\phi_1)}{\phi_{11}'} \ \mbox{and}\
j^k \epsilon_2 =-j^k \frac{\eta_2(\phi_2)}{\phi_{22}'}$ if $\delta=1$ or $j^k \epsilon_2=-j^k \frac{\eta_1(\phi_2)}{\phi_{21}'}$ if $\delta=2$. So,
$$h_1=\frac{\eta_2(\phi_1)\phi_{11}'-\eta_1(\phi_1)\phi_{12}'}{m_1t_1^{m_1-1}}\ mod\ t_1^{k+1}\ \
\mbox{and}\ \
h_2=\frac{\eta_2(\phi_2)\phi_{21}'-\eta_1(\phi_2)\phi_{22}'}{m_2t_2^{m_2-1}} \ mod\ t_2^{k+1}.$$ Defining 
$\omega=\eta_2dx-\eta_1dy\in\Omega[\delta]$, we have that $h_i=j^k \frac{\omega\circ\phi_i}{m_it_i^{m_i-1}}$, $i=1,2$.

Conversely, given $\omega=g_2dx+g_1dy\in\Omega[\delta]$ where
$g_1=\beta(\delta-1)y+h$ with $h,g_2\in (x,y)^2\mathbb{C}\{x,y\}$ and
$\beta\in\mathbb{C}$, consider
$\eta_1=-g_1,\eta_2=g_2,\epsilon_1=\frac{g_1(\phi_1)}{m_1t_1^{m_1-1}}\in (t_1)^2\mathbb{C}\{t_1\},\epsilon_2=-\frac{g_2(\phi_2)}{m_2t_2^{m_2-1}}$ if $\delta=1$ or $\epsilon_2=\frac{g_1(\phi_2)}{m_2t_2^{m_2-1}}$ \ if $\delta=2$. So, from (\ref{tangentvector}), we have that 
$[(0,h_1),(-1)^{\delta}((2-\delta)h_2,(\delta-1)h_2)]\in
T\A^k[\delta] (\phi)$, where
$h_i=j^k \frac{\omega(\phi_i)}{m_it_i^{m_i-1}}$, $i=1,2$. \cqd

In the sequel we will need the notions of fibers $F_i$ and the set $M$ of maximal points of the sets $\Lambda[\delta]$, which are defined in a similar way as for $\Lambda$.  We will also use the notation $\underline{k}= (k,k)\in\mathbb{N}^2$.

\begin{corollary}
\label{Diag} Let $\phi=[\phi_1,\phi_2] \in
\Sigma_{\Gamma,\Lambda[\delta]}$ and $k \in
\mathbb{N}$.
\begin{enumerate}
\item[\rm (a)] If $k>m_1$ then $F_1(\underline{k}) \neq \emptyset$ if and only if  $[(0,t_1^{k}),(0,0)] \in T\A^k[\delta](\phi)$;
\item[\rm (b)] If $k>m_2$ then $F_2(\underline{k}) \neq \emptyset$ if and only if  $[(0,0),(-1)^{\delta}((2-\delta)t_2^{k},(\delta-1)t_2^k)] \in T\A^k[\delta](\phi)$;
\item[\rm (c)] If $\underline{k} \in M$ then there exist $a,b \in
\C^*,$ such that  $$[(0,at_1^{k}),(-1)^{\delta}((2-\delta)bt_2^k,(\delta-1)bt_2^k)]\in T\A^k[\delta](\phi).$$
\end{enumerate}
\end{corollary}
\pf We have that $(\gamma_1,\gamma_2)\in\Lambda[\delta]$ if and only if there exists  $\omega\in\Omega[\delta]$
such that $ord_{t_i}\frac{\omega(\phi_i)}{m_it_i^{m_i-1}}=\gamma_i$. This, in turn, is equivalent, from the preceding result, to
\begin{equation}
[( 0,h_1),(-1)^{\delta}((2-\delta)h_2,(\delta-1)h_2)] \in T\A^k[\delta](\phi), \label{ajuda}
\end{equation}
where $h_i=j^k \frac{\omega(\phi_i)}{m_it_i^{m_i-1}}$.

Now, suppose that $k>m_1$. Then $F_1(\underline{k})\neq\emptyset$ if and only if there exists
$(k,\gamma)\in\Lambda[\delta]$ with $\gamma>k$. The last condition, from 
(\ref{ajuda}), is equivalent to the condition $[(0,t_1^{k}),(0,0)] \in T\A^k[\delta](\phi)$, proving in this way $(a)$. The proof of $(b)$ is analogous.

Now, if $\underline{k}\in M$, then from (\ref{ajuda}) we have that 
$[(0,at_1^{k}),(-1)^{\delta}((2-\delta)bt_2^k,(\delta-1)bt_2^k)]\in T\A^k[\delta](\phi)$.

\cqd

The next result will give us the normal forms of bigerms under the action of the  group $\A[\delta]$.

\begin{proposition}
Let $\phi=[\phi_1,\phi_2] \in \Sigma_{\Gamma,\Lambda[\delta]}$. If 
$F_i(\underline{k}) \neq \emptyset$ and $k>m_i$ for some $i\in\{1,2\}$ (respectively $\underline{k}\in M$), then there exists $\varphi=[\varphi_1, \varphi_2]\in \Sigma_{\Gamma,\Lambda[\delta]}$ such that 
$\varphi$ is $\A[\delta]$-equivalent to $\phi$ with $j^{k-1}\varphi =
j^{k-1}\phi$ and $j^{k}\varphi_i=j^{k-1}\phi_i$ (respectively $j^k\varphi_1=j^{k-1}\phi_1$ or $j^k\varphi_2=j^{k-1}\phi_2$).
\label{normal}
\end{proposition}
\pf From Corollary \ref{Diag}(a), if $F_1(\underline{k})
\neq \emptyset$ and $k > m_1$, then $[(0,t_1^{k}),(0,0)] \in H_\phi^k[\delta]\cap T\A^k[\delta](\phi)$. It follows from CTT that $j^{k-1}\phi$ is $\A^k[\delta]$-equivalent to $[j^k\phi_1,j^{k-1}\phi_2]$ and therefore there exists $\varphi$ which is $\A[\delta]$-equivalent to $\phi$ such that $j^k\varphi_1=j^{k-1}\phi_1$ and 
$j^{k-1}\varphi= j^{k-1}\phi$. The case
$F_2(\underline{k}) \neq \emptyset$ and $k > m_2$ is analogous.

If $\underline{k}\in M$, then $F_1(\underline{k})=F_2(\underline{k})= \emptyset$ and, from  Corollary \ref{Diag}(c), the element given by
$[(0,dat_1^{k}),(-1)^{\delta}((2-\delta)dbt_2^k,(\delta-1)bdt_2^k)]$ with well determined $a,b \in \C^*$ and arbitrary  $d\in\mathbb{C}$ belongs to $H_\phi^k[\delta]\cap T\A^k[\delta](\phi)$. Choosing $d$ conveniently, it follows, as we argued before, that there exists a bigerm $\varphi$ which is $\A[\delta]$-equivalent to $\phi$ such that 
$j^{k-1}\varphi = j^{k-1}\phi$ with 
$j^k\varphi_1=j^{k-1}\phi_1$ or $j^k\varphi_2=j^{k-1}\phi_2$ according to the choice of $d$.
\cqd

Since there are two different choices to be made in this process when $\underline{k}$ is in $M$, given  $\phi=[\phi_1,\phi_2]\in\Sigma_{\Gamma,\Lambda[\delta]}$ and $\underline{k}\in M$ we will choose an $\A[\delta]$-equivalent $\varphi$  to $\phi$ such that 
$j^{k-1}\varphi = j^{k-1}\phi$ and
$j^k\varphi_1=j^{k-1}\phi_1$. In this way, we have the following description of the normal forms for bigerms in $\Sigma_{\Gamma,\Lambda[\delta]}$:

\begin{theorem}{\sc ($\A[\delta]$-normal form)}
A bigerm $\phi=[\phi_1,\phi_2] \in
\Sigma_{\Gamma,\Lambda[\delta]}$ is always $\A[\delta]$-equivalent to  a $\varphi=[\varphi_1,\varphi_2]$ such that 
\begin{equation}
\label{FormaNormal}
\begin{array}{ccc}
\varphi_1 =\left (t_1^{m_1},\displaystyle \sum_{\underline{j}\not\in M\atop
F_1(\underline{j})=\emptyset} a_{1j}t_1^{j}\right ) & &
\varphi_2 =\left \{
\begin{array}{lc}
\left (\displaystyle\sum_{F_2(\underline{j})=\emptyset}
a_{2j}t_2^{j},t_2^{m_2} \right ) & \mbox{if}\ \delta=1 \\ \\
\left (t_2^{m_2},\displaystyle\sum_{F_2(\underline{j})=\emptyset}
a_{2j}t_2^{j}\right ) & \mbox{if}\ \delta=2\;.
\end{array}
\right .
\end{array}
\end{equation}
\end{theorem}

Now, we will prove the uniqueness of the $\A[\delta]$-normal form, by arguments similar to those used in \cite{HH}.

The set 
$$N=\{\varphi \in\Sigma_{\Gamma ,\Lambda[\delta]}; \ \varphi \ \mbox{as given in} \ (\ref{FormaNormal})\}$$
is an open set in some affine space of finite dimension. Denoting by $N^k$ the space $j^k(N)$, we have the following 
lemma:

\begin{lemma} If $\phi=[\phi_1,\phi_2] \in N$, then for all $k>min\{m_1,m_2\}$, we have
$$N^k \cap \{j^k\phi +T\A^k[\delta](j^k \phi) \}=\{j^k\phi\}.$$
\end{lemma}
\pf Suppose the assertion not true. Take $k$ minimal with the
following property: 
$$N^k \cap \{j^k\phi +T\A^k[\delta](j^k \phi) \}\neq\{j^k\phi\}.$$

So, there exists $\psi \in N^k \cap \{j^k\phi + T\A^k[\delta](j^k\phi)\}$ such that $\psi
\neq j^k\phi$ and $j^{k-1}\psi=j^{k-1}\phi$ because
$k$ is minimal. Therefore, there exist $b_1,b_2 \in {\mathbb C}$
with $b_1\neq 0$ or $b_2\neq 0$ such that
\[
\psi-j^k\phi=[(0,b_1t_1^k),((2-\delta)b_2t_2^k,(\delta-1)b_2t_2^k)]\in T\A^k[\delta](j^k\phi).
\]

If $F_i(\underline{k})\neq \emptyset$ for some $i=1,2$, then we have
a contradiction, since $\psi, j^k\phi \in N^k$ are given as in (\ref{FormaNormal}). So, we have
$\underline{k} \in M$. But since $\psi, j^k\phi \in N^k$ we
have $b_1 =0$, then $b_2 \neq 0$. In this way,     $\psi-j^k\phi=[(0,0),((2-\delta)b_2t_2^k,(\delta-1)b_2t_2^k)]
\in T\A^k[\delta](j^k\phi)$, and
$F_2(\underline{k}) \neq \emptyset$ which is again a contradiction. \cqd

Now we conclude the proof of the uniqueness of the $\A[\delta]$-normal forms.

Let $\phi=[\phi_1,\phi_2] \in \Sigma_{\Gamma,\Lambda[\delta]}$.
Observe that our bigerms are finitely determined up to order $c$ (as defined at the beginning of this section), that is, $\phi$ is $\A[\delta]$-equivalent to $j^c \phi$. We have to prove that $N^{c}\cap \A^{c}[\delta](j^{c}\phi) = \{j^{c}\phi\}$. Suppose that $\varphi\in
N^{c}\cap\A^{c}[\delta](j^{c}\phi)$, with $\varphi \neq
j^{c}\phi$. Since $\A^{c}[\delta](j^{c}\phi)$ is arcwise
connected, there exists an arc in $\A^{c}[\delta](j^{c}\phi)$
joining $j^{c}\phi$ to $\varphi$. Since the reduction process to the normal
form is continuous, it follows that $j^{c}\phi$ is not an
isolated point in $N^{c}\cap\A^{c}[\delta](j^{c}\phi)$.
This is a contradiction because of Lemma 7.\cqd

Since the $\A$-action on bigerms is  the composition of the $\A[\delta]$-action with homotheties, the $\A$-action on the $\A[\delta]$-normal forms reduces to the action of the group of homotheties. 

\section{Homothety Action}

We will consider initially the case of bigerms with transversal components. 

In this case, we may write
$$\phi=[\phi_1,\phi_2]=\left[\left(t_1^{m_1}, \sum_{j=j_1}^{c} a_{1j}t_1^j \right), \ \ \left(
\ \sum_{j=j_2}^{c} a_{2j}t_2^j, t_2^{m_2}\right)\right].$$

In order to preserve the above form, we have to consider the following particular homotheties:
 $(\rho_1,\rho_2,\sigma)\in\mathcal H$ with $\sigma
(x,y)=(\alpha_1 x,\alpha_2 y)$ and $\rho_i
(t_i)=\alpha_i^{\frac{1}{m_i}}t_i$,  $\alpha_i\in \C^*$, $i=1,2$. 

In this way get 
$$(\rho_1,\rho_2,\sigma)\cdot \phi = \left[ \left(t_1^{m_1}, \sum_{j=j_1}^{c} \alpha_1^{-\frac{j}{m_1}}\alpha_2a_{1j}t_1^j \right), \; \left(\sum_{j=j_2}^{c} \alpha_2^{-\frac{j}{m_2}}\alpha_1a_{2j}t_2^j, \ t_2^{m_2}
\right)\right].$$

In this situation, with a convenient choice of $\alpha_1$ and $\alpha_2$ we may reduce two any non-zero coefficients in the above sums to $1$. We will always choose to apply this reduction to the coefficients of the terms lower order of $\phi_1$, if they exist. If not, we continue in the same way the reduction on the terms of $\phi_2$.

Similarly, when the components of $\phi$ have same tangent, that is, when 
$$\phi=[\phi_1,\phi_2]=\left[\left(t_1^{m_1}, \sum_{j=j_1}^{c} a_{1j}t_1^j \right), \ \ \left(t_2^{m_2},
\ \sum_{j=j_2}^{c} a_{2j}t_2^j \right)\right],$$
we have to consider $\sigma (x,y)=(\alpha_1 x,\alpha_2 y)$ and
$\rho_i(t_i)=\alpha_1^{\frac{1}{m_i}}t_i$ with 
$\alpha_i\in\mathbb{C}*$, $i=1,2$. In this case, we get
$$\sigma\circ\phi_i\circ\rho_i^{-1}(t_i)=\left(t_i^{m_i}, \sum_{j=j_i}^{c} \alpha_1^{-\frac{j}{m_{i}}}
\alpha_2a_{ij}t_i^j \right), i=1,2.$$

In the same way as above, we may reduce to $1$ any two coefficients in the above sums, unless both components of $\phi_1$ and $\phi_2$ are monomials with $m_1=m_2$ and $j_1=j_2$. In this case, we may reduce to $1$ only one of the coefficients.

The above discussion may be summarized in the following theorem:

\begin{theorem}\label{normalforms}
Any $\phi \in \Sigma_{\Gamma,\Lambda[\delta]}$ is $\A$-equivalent to one in the following form:

\noindent{\bf Distinct tangents case}

\begin{enumerate}
\item[\rm a)] $\left[\left(t_1^{m_1}, t_1^{j_1}+t_1^k+ \displaystyle\sum_{\underline{j}\not\in M\atop F_1(\underline{j})
=\emptyset} a_{1j}t_1^j \right), \ \
\left(\displaystyle\sum_{F_2(\underline{j}) = \emptyset} a_{2j}t_2^j
\ , \ t_2^{m_2} \ \right)\right];$

\item[\rm b)] $\left[\left(t_1^{m_1}, t_1^{j_1}\right),
(t_2^{j_2}, t_2^{m_2})\right];$

\item[\rm c)] $\left[(t_1^{m_1}, t_1^{j_1}),(0, t_2)\right];$

\item[\rm d)] $\left[ (t_1,0), (0,t_2)\right]$.

\end{enumerate}

\noindent{\bf Same tangents case}
\begin{enumerate}
\item[\rm a$'$)] $\left[\left(t_1^{m_1},
t_1^{j_1}+t_1^k+\displaystyle\sum_{\underline{j}\not\in M\atop
F_1(\underline{j})=\emptyset} a_{1j}t_1^j \right), \ \
\left(t_2^{m_2}, \ \displaystyle\sum_{F_2(\underline{j})= \emptyset}
a_{2j}t_2^j \right)\right];$

\item[\rm b$'$)] $\left[ (t_1^{m_1}, t_1^{j_1}), (t_2^{m_2},
t_2^{j_2})\right]$ with $m_1 \neq m_2 \ \mbox{or}\ j_1 \neq j_2$;

\item[\rm c$'$)] $\left[ (t_1^{m_1}, t_1^{j_1}), (t_2^{m_1},
at_2^{j_1})\right]$, with $a\not\in \{0,1\} $.

\item[\rm d$'$)] $[(t_1^{m_1},t_1^{j_1}),(t_2,0)].$
\end{enumerate}
\end{theorem}

Let us remark that two bigerms in the above list with distinct normal forms are not $\A$ equivalent since their corresponding sets $\Lambda$ are not equal.

In what follows we will describe the homotheties that preserve the above normal forms. Since in cases b), c), d), b$'$), c$'$) and d$'$), the homotheties act as the identity, we have only to describe such homotheties in the remaining cases a) and a$'$). In these cases $\sigma (x,y)=(\alpha^{m_1}x,\alpha^{j_1}y)$, $\rho_1(t_1)=\alpha t_1$, with $\alpha^{k-j_1}=1$ and

\vspace{0.2cm} \noindent{\bf Case a)}: $\rho_2(t_2)=\alpha^{\frac{j_1}{m_2}}t_2$. In this case, two bigerms with coefficients $a_{ij}$ and $b_{ij}$ are $\H$-equivalent if and only if 
\[
a_{1j}\alpha^{j_1-j}=b_{1j}, \ \text{and}  \  a_{2j}\alpha^{\frac{m_1m_2-j_1j}{m_2}}=b_{2j}.
\]

\vspace{0.2cm} \noindent{\bf Case a$'$)}: $\rho_2(t_2)=\alpha^{\frac{m_1}{m_2}}t_2$. In this case, two bigerms with coefficients $a_{ij}$ and $b_{ij}$ are $\H$-equivalent if and only if 
\[
 a_{ij}\alpha^{\frac{j_1m_i-j m_1}{m_i}}=b_{ij}, \ \ \ i=1,2.
\]

\section{Final Remarks}

Given any two bigerms $\phi=[\phi_1,\phi_2]$ and
$\psi=[\psi_1,\psi_2]$, to verify if they are $\A$-equivalent we may proceed as follows:
\begin{enumerate}
\item If semigroup $\Gamma_\phi$ and the semigroups $\Gamma^1_\psi$ and $\Gamma^2_\psi$ corresponding to the two possible orders of the branches of $\psi$ are such that $\Gamma^1_{\psi} \neq \Gamma_\phi \neq \Gamma^2_{\psi}$, then $\phi$ and $\psi$ are not $\A$-equivalent. If this is not the case, choose the order of the branches of $\psi$ to force the equality of the semigroups of $\phi$ and $\psi$.  

\item If $\Lambda_{\phi}[\delta]\neq\Lambda_{\psi}[\delta]$, then the bigerms are not  $\A$-equivalent.

\item If $\Lambda_{\phi}[\delta]=\Lambda_{\psi}[\delta]$, we take representatives for $\phi$ and $\psi$ in normal form as in Theorem 8. \label{repete1}

\item We verify if one of the homotheties that preserve the normal form transforms $\phi$ into $\psi$. In such case, the two bigerms are $\A$-equivalent.\label{repete2}

\item If this is not the case, we have to permutate the branches of one of the bigerms and repeat steps \ref{repete1} and \ref{repete2}.
\end{enumerate}

To give an explicit example for a pair of bigerms as in step 5 above, consider 
$$\phi=[(t_1^{m},t_1^{j}),(t_2^{m},at_2^{j})] \ \ \text{and}  \ \ 
\psi=[\psi_1,\psi_2]=[(t_1^{m},t_1^{j}),(t_2^{m},bt_2^{j})], \ \ a,b\not \in \{0,1\},$$ which are in normal form c$'$). 

So, $\phi$ and $\psi$ are $\H$-equivalent if and only if $a=b$. On the other hand, if we permutate the branches of $\psi$ and put it in normal form, we get 
$[\psi_2,\psi_1]=[(t_2^{m},t_2^{j}),(t_1^{m},\frac{1}{b}t_1^{j})]$. Therefore, $\phi$ and $\psi$ are $\A$-equivalent if and only if $a=b$ or $a=\frac{1}{b}$.  
This is a generalization of Example 3 of \cite{CDF}.

In the irreducible case, in each equisingularity class determined by semigroups of the form $\mathbb N$, $\langle 2,j \rangle$ with $j\equiv 1 \bmod 2$, or $\langle 3, 3+\alpha\rangle$
with $\alpha=1,2$, all curves are analytically equivalent to a monomial curve, that is, for any of these equisingularity classes we have one possible set $\Lambda$, namely, $\Lambda=\Gamma \setminus \{0\}$.

Using the description of the semigroup, the set of maximal points as described in \cite{Ga} and doing some computations with differentials, we get the following table for bigerms with  transversal components and whose semigroups are as described above.

\bigskip

\begin{tabular}{|l|c|}
\hline {\bf $(m_1,m_2)$} & Normal Form \\
\hline $(1,1)$ &  $(t_1,0)\ (0,t_2)$ \\ 
\hline $(1,2)$ &
$(t_1^2,t_1^{j})\ (0,t_2)$ \ \ $j\equiv1 \bmod 2$\\
\hline $(1,3)$ & $(t_1^3,t_1^{3+\alpha})\ (0,t_2);\ \alpha=1,2$ \\
  & $(t_1^3,t_1^{3+\alpha}+t_1^{3+2\alpha})\ (0,t_2);\ \alpha=1,2$ \\
\hline $(2,2)$ & $(t_1^2,t_1^{j_1})\ (t_2^{j_2},t_2^2)$  \  \ $j_i\equiv 1 \bmod 3$, \ $i=1,2$\\
\hline $(2,3)$
& $(t_1^3,t_1^{3+\alpha})\ (t_2^{j},t_2^2);\ \alpha=1,2$, \ $j\equiv 1 \bmod 2$ \\
  & $(t_1^3,t_1^{3+\alpha}+t_1^{3+2\alpha})\ (at_2^{j},t_2^2);\ \alpha=1,2,\ a\neq 0$, \ $j\equiv 1 \bmod 2$ \\
\hline   & $(t_1^3,t_1^{3+\alpha_1})\ (t_2^{3+\alpha_2},t_2^3);\ \alpha_1,\alpha_2=1,2$ \\
 $(3,3)$ & $(t_1^3,t_1^{3+\alpha_1}+t_1^{3+2\alpha_1})\ (at_2^{3+\alpha_2},t_2^3);\ \alpha_1,\alpha_2=1,2,\ a\neq 0$ \\
 & $(t_1^3,t_1^{3+\alpha_1}+t_1^{3+2\alpha_1})\ (at_2^{3+\alpha_2}+bt_2^{3+2\alpha_2},t_2^3);\ \alpha_1,\alpha_2=1,2,\ a,b\neq 0$ \\
\hline

\end{tabular}

\bigskip

\noindent Authors Affiliations:\medskip

\noindent Abramo Hefez (hefez@mat.uff.br)\medskip

\noindent {\bf Universidade Federal Fluminense}\bigskip

\noindent Marcelo Escudeiro Hernandes (mehernandes@uem.br) 
and \medskip

\noindent Maria Elenice Rodrigues Hernandes (merhernandes@uem.br)\medskip

\noindent {\bf Universidade Estadual de Maring\'a}
\end{document}